# COMPETITION INTERFACES AND SECOND CLASS PARTICLES[1]

By Pablo A. Ferrari and Leandro P. R. Pimentel

*Universidade de São Paulo*

The one-dimensional nearest-neighbor totally asymmetric simple exclusion process can be constructed in the same space as a last-passage percolation model in $\mathbb{Z}^2$. We show that the trajectory of a second class particle in the exclusion process can be linearly mapped into the competition interface between two growing clusters in the last-passage percolation model. Using technology built up for geodesics in percolation, we show that the competition interface converges almost surely to an asymptotic random direction. As a consequence we get a new proof for the strong law of large numbers for the second class particle in the rarefaction fan and describe the distribution of the asymptotic angle of the competition interface.

**1. Introduction.** The relation between the totally asymmetric nearest-neighbor simple exclusion process in dimension one and two-dimensional last-passage percolation models is well known since the seminal work of Rost [19]. The macroscopic behavior of the density profile of the exclusion process is governed by the Burgers equation [1, 17]. This corresponds to the "shape theorem" in last-passage percolation [19, 20]. An important property of the exclusion process is that the so-called second class particles (that follow roughly the behavior of a perturbation of the system) are asymptotically governed by the characteristics of the Burgers equation. When there is only one characteristic, the second class particle follows it [5, 18, 21]; when there are infinitely many, the particle chooses one of them at random to follow [6]. These results hold when the initial distribution is a product measure with densities $\lambda \in (0,1]$, $\rho \in [0,1)$, to the left and right of the origin, respectively. The existence of infinitely many characteristics occurs at points where the solution of the Burgers equation is a rarefaction front. The rescaled position of the second class particle converges almost surely to a random variable

Received June 2004; revised September 2004.
[1]Supported in part by CNPq, FAPESP and PRONEX.
*[AMS 2000 subject classifications.](AMS 2000 subject classifications.)* Primary 60K35; secondary 82B.

uniformly distributed in the interval $[1-2\lambda, 1-2\rho]$ as time goes to infinity [14]. A similar phenomenon has been observed in first-passage percolation starting from two growing clusters competing for space: the rescaled competition interface converges almost surely to a random direction [16] with a so far unknown distribution. Motivated by this we investigate the relation between the second class particle and the competition interface in last-passage percolation. We conclude that one object can be mapped into the other (as processes) realization by realization. Indeed, the difference of the coordinates of the competition interface at time $t$ is exactly the position of the second class particle at that time (see Proposition 3 and Lemma 6). We show a law of large numbers for the competition interface in the positive quadrant $(\mathbb{Z}^+)^2$; this corresponds to $\lambda = 1$ and $\rho = 0$. Our mapping then permits to describe the distribution of the angle of the competition interface in last-passage percolation (Theorem 1) and to give a new proof of the strong law of large numbers for the second class particle (Theorem 2, for the moment restricted to the case $\lambda = 1$ and $\rho = 0$; we comment in the final remarks what should be done in the other cases). A key tool to prove the asymptotic behavior of the competition interface is the study of the *geodesics*, random paths maximizing the passage time. We show that each semi-infinite geodesic has an asymptotic direction and that two semi-infinite geodesics with the same direction must coalesce. The approach follows Newman [15] who proved analogous results for first-passage percolation (see also [9, 10]).

In Section 2 we introduce the models, state the results and prove them. In Section 3 we show properties of the geodesics needed for the proofs.

**2. Last-passage percolation and simple exclusion.** Let $\mathcal{W} = (w(z), z \in \mathbb{Z}^2)$ be a family of independent random variables with exponential distribution of mean 1. Let $\mathbb{P}$ and $\mathbb{E}$ be the probability and expectation induced by these variables in the product space $\Omega = (\mathbb{R}^+)^{\mathbb{Z}^2}$.

Given $z = (i, j)$, $z' = (i', j')$ in $\mathbb{Z}^2$ with $i \leq i'$ and $j \leq j'$, we say that $(z_k, k = 1, \ldots, n)$ is an *up/right path from $z$ to $z'$* if $z_1 = z$, $z_n = z'$ and $z_{k+1} - z_k \in \{(0,1), (1,0)\}$ for $k = 1, \ldots, n-1$. Let $\Pi(z, z')$ be the set of up/right paths from $z$ to $z'$. The *maximal length* between $z$ and $z'$ is defined by

$$(1) \qquad G(z, z') := \max_{\pi \in \Pi(z, z')} \left\{ \sum_{z'' \in \pi} w(z'') \right\}.$$

This model is called *last-passage percolation*. Since we are interested in the paths starting at $(1,1)$, we use the notation $G(z) = G((1,1), z)$. This function satisfies the recurrence relation

$$(2) \qquad G(z) = w(z) + \max\{G(z - (0,1)), G(z - (1,0))\}$$

with $G(i,j) = 0$ if either $i = 0$ or $j = 0$. We say that a point $z$ is *infected* at time $t$ if $z \in \mathbf{G}_t$, where

$$\mathbf{G}_t := \{z \in (\mathbb{Z}^+)^2 : G(z) \leq t\}$$

is called the *infected region*. Let $Q(i,j) := (i-1, i] \times (j-1, j]$ be the unit square having $(i,j)$ as north-east vertex. The set $\overline{\mathbf{G}}_t := \bigcup_{z \in \mathbf{G}_t} Q(z)$ describes the subset of $(\mathbb{R}^+)^2$ attained by the infection at time $t$. The random process $\overline{\mathbf{G}}_t$ is called a *spatial growth model* and describes a growing Young tableau. The *growth interface* is defined by

(3) $\qquad \gamma_t := \{(i,j) \in (\mathbb{Z}^+)^2 : G(i,j) \leq t \text{ and } G(i+1, j+1) > t\}.$

The polygonal curve interpolating the points of $\gamma_t$ that are at distance 1 separates the infected region $\overline{\mathbf{G}}_t$ and its complement.

Rost [19] proved a "shape theorem" for $\overline{\mathbf{G}}_t$: with $\mathbb{P}$ probability 1, for all $\varepsilon > 0$ there exists a $t_0$ such that for all $t > t_0$,

(4) $\qquad\qquad\qquad t(1-\varepsilon)\mathcal{M} \subset \overline{\mathbf{G}}_t \subset t(1+\varepsilon)\mathcal{M}$

where $\mathcal{M} := \{(u,v) \in (\mathbb{R}^+)^2 : \mu(u,v) \leq 1\}$ and

(5) $\qquad\qquad\qquad \mu(u,v) := (\sqrt{u} + \sqrt{v})^2.$

The interface $\gamma_t$ converges to $\{(u,v) : \mu(u,v) = 1\}$ in the same sense: with $\mathbb{P}$ probability 1, for all $\varepsilon > 0$ there exists a $t_0$ such that for all $t > t_0$,

(6) $\qquad\qquad\qquad \gamma_t \subset [t(1+\varepsilon)\mathcal{M}] \setminus [t(1-\varepsilon)\mathcal{M}].$

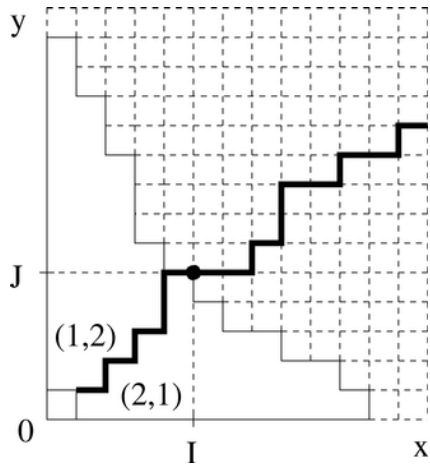

Fig. 1. *Growth and competition interfaces.*

*Competing spatial growth.* The sets of points infected through $(2,1)$ and $(1,2)$, respectively, are defined by

$$\mathbf{G}_t^{21} := \{z \in (\mathbb{Z}^+)^2 : G(z) \leq t \text{ and } G(z) = w(1,1) + G((2,1),z)\},$$

$$\mathbf{G}_t^{12} := \{z \in (\mathbb{Z}^+)^2 : G(z) \leq t \text{ and } G(z) = w(1,1) + G((1,2),z)\}.$$

The process $(\mathbf{G}_t^{21}, \mathbf{G}_t^{12})$ describes a competing spatial growth model between two different infections (see Figure 1). For related models in first-passage percolation see [3, 8, 16]. One can see that the regions $\mathbf{G}_t^{21}, \mathbf{G}_t^{12}$ are connected, $\mathbf{G}_t = \{(1,1)\} \cup \mathbf{G}_t^{21} \cup \mathbf{G}_t^{12}$ and that the *competition interface* $\varphi = (\varphi_0, \varphi_1, \ldots)$ between $\mathbf{G}_\infty^{21}$ and $\mathbf{G}_\infty^{12}$ can be defined inductively as follows: $\varphi_0 = (1,1)$ and for $n \geq 0$,

(7) $$\varphi_{n+1} = \begin{cases} \varphi_n + (1,0), & \text{if } \varphi_n + (1,1) \in \mathbf{G}_\infty^{21}, \\ \varphi_n + (0,1), & \text{if } \varphi_n + (1,1) \in \mathbf{G}_\infty^{12}. \end{cases}$$

So that, if we paint blue the squares $Q(z)$ with $z \in \mathbf{G}_\infty^{21}$ and red the squares $Q(z)$ with $z \in \mathbf{G}_\infty^{12}$, the line obtained by linear interpolation of $\varphi_0, \varphi_1, \ldots$ separates the blue and red regions. The square $Q(1,1)$ gets no color. Definition (7) is equivalent to

(8) $$\varphi_{n+1} = \arg\min\{G(\varphi_n + (1,0)), G(\varphi_n + (0,1))\}, \qquad n \geq 0.$$

Note that given $G(z)$ for all $z$, the interface $\varphi$ chooses locally the shorter step to go up or right. We prove that $\varphi$ has an asymptotic (random) direction and compute the law of the direction:

THEOREM 1.

(9) $$\lim_{n \to \infty} \frac{\varphi_n}{|\varphi_n|} = e^{i\theta}, \qquad \mathbb{P}\text{-}a.s.$$

where $\theta = \theta(\mathcal{W})$ is a random angle in $[0, 90°]$ with law

(10) $$\mathbb{P}(\theta \leq \alpha) = \frac{\sqrt{\sin \alpha}}{\sqrt{\sin \alpha} + \sqrt{\cos \alpha}}.$$

*Second class particles in simple exclusion.* The one-dimensional nearest-neighbor totally asymmetric simple exclusion process is a Markov process $(\eta_t, t \geq 0)$ in the state space $\{0,1\}^\mathbb{Z}$. $\eta_t(x)$ indicates if there is a particle at site $x$ at time $t$; only one particle is allowed at each site. At rate 1, if there is a particle at site $x \in \mathbb{Z}$, it attempts to jump to $x+1$; if there is no particle in $x+1$ the jump occurs, otherwise nothing happens. To construct a realization of this process à la Harris, one considers independent one-dimensional Poisson processes $\mathcal{N} = (N_x(\cdot), x \in \mathbb{Z})$ of intensity 1; let $\mathbb{Q}$ be the law of $\mathcal{N}$. The process $(\eta_t, t \geq 0)$ can be constructed as a deterministic function of the initial configuration $\eta$ and the Poisson processes $\mathcal{N}$ as follows: if $s$ is a Poisson epoch of $N_x$ and there is a particle at $x$ and no particle at $x+1$ in the configuration $\eta_{s-}$, then at time $s$ the new configuration is obtained by making the particle jump from $x$ to $x+1$. This construction is well defined; see [4], for instance. Let $\Phi$ be the function that takes $\eta$ and $\mathcal{N}$ to $(\eta_t,\ t\geq 0)$. Let $\eta^0$ and $\eta^1$ be two arbitrary configurations. The *basic* coupling between two exclusion processes with initial configurations $\eta^0$ and $\eta^1$ respectively is the joint realization $(\Phi(\eta^0,\mathcal{N}),\Phi(\eta^1,\mathcal{N})) = ((\eta_t^0,\eta_t^1),\ t\geq 0)$ obtained by using the same Poisson epochs for the two different initial conditions. Liggett [11, 12] are the default references for the exclusion process.

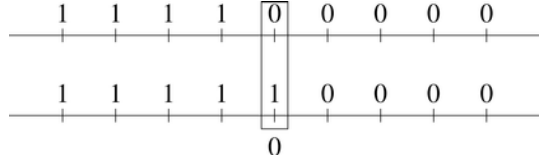

Fig. 2. *Second class particle. The first line is $\eta^0$ and the second one is $\eta^1$.*

Let $\eta^0$ and $\eta^1$ be two configurations defined by

(11) $$\eta^0(x) = \mathbf{1}\{x \leq -1\}, \qquad \eta^1(x) = \mathbf{1}\{x \leq 0\}.$$

These configurations are full to the left of the origin and empty to the right of it and differ only at the origin (see Figure 2). Call $X(0) = 0$ the site where both configurations differ at time zero. With the basic coupling, the configurations at time $t$ differ only at the site $X(t)$ defined by

$$X(t) := \sum_x x\mathbf{1}\{\eta_t^0(x) \neq \eta_t^1(x)\}.$$

$(X(t),\ t \geq 0)$ is the trajectory of a "second class particle." The process $((\eta_t^0, X(t)), t \geq 0)$ is Markovian but the process $(X(t), t \geq 0)$ is not. The motion of $X(t)$ depends on the configuration of $\eta_t^1$ in its neighboring sites. The second class particle jumps one unit to the right at rate 1 if there is no $\eta^1$ particle in its right nearest neighbor and it jumps one unit to the left at rate 1 if there is an $\eta^1$ particle in its left nearest-neighbor site, interchanging positions with it. Ferrari and Kipnis [6] proved that $X(t)/t$ converges in distribution to a uniform random variable as $t \to \infty$ for initial configurations distributed according to product measures with densities $\lambda > \rho$ to the left and right of the origin, respectively. In these cases, Mountford and Guiol [14] proved almost sure convergence. Our approach gives an alternative proof to Mountford and Guiol in the case $\lambda = 1$ and $\rho = 0$:

THEOREM 2. *Let $(X(t),\ t \geq 0)$ be the trajectory of a second class particle put initially at the origin in the one-dimensional totally asymmetric nearest-neighbor simple exclusion process starting with the configuration $\eta^1$ defined* by $\eta^1(x) = \mathbf{1}\{x \leq 0\}$. Then

$$\lim_{t \to \infty} \frac{X(t)}{t} = U, \qquad \mathbb{Q}\text{-}a.s. \tag{12}$$

where $U = U(\mathcal{N})$ is a random variable with uniform distribution in $[-1, 1]$.

*Pair representation of the second class particle.* It is convenient to represent the second class particle with a pair hole–particle. For that we consider the initial configuration $\eta^{01}$ defined by

$$\eta^{01}(x) = \begin{cases} \eta^1(x), & \text{if } x \leq -1, \\ \eta^1(x-1), & \text{if } x > 1, \\ 0, & \text{if } x = 0, \\ 1, & \text{if } x = 1. \end{cases} \tag{13}$$

This configuration has a particle at site 1 called *particle and a hole at site 0 called *hole. The pair *hole–*particle is called *pair (see the configuration before jump in Figure 3). The process $\eta_t^{01}$ is constructed using the Poisson marks as before; ignoring the *pair, the process is just the exclusion process starting with the configuration $\eta^{01}$. On top of it we define the evolution of the *pair as follows: when a particle (from the left) jumps over the *hole, the *pair moves one unit to the left (giving rise to the configuration after the jump in Figure 3); when the *particle jumps to the right (over a hole), the *pair moves one unit to the right. This is the same behavior as that of the second class particle; the difference is that the second class particle occupies only one site while the *pair occupies two sites. Call $P^*(t)$ and $H^*(t)$ the position of the *particle and *hole respectively at time $t$; clearly $P^*(t) = H^*(t) + 1$ for all times. If we collapse again the *pair to one site by defining $\bar{\eta}_t(x) = \eta_t^{01}(x)$ for $x < H^*(t)$, $\bar{\eta}_t(x) = \eta_t^{01}(x+1)$ for $x \geq P^*(t)$, then

(14) the process $(\bar{\eta}_t, H^*(t),\ t \geq 0)$ has the same law as $(\eta_t^1, X(t),\ t \geq 0)$.

In Lemma 6 we give an explicit construction which maps these processes for almost all realizations.

*Growth model and simple exclusion.* Rost [19] showed that the simple exclusion process can be constructed in the probability space induced by $\mathcal{W}$, where the oriented percolation model is defined. This can be done for any initial configuration; we do it for the process with initial configuration $\eta^{01}$ as follows. Let

$$\begin{aligned} P_1(0) &= 1, \qquad H_1(0) = 0; \\ P_i(0) &= -i+1 \quad \text{and} \quad H_i(0) = i, \qquad i \geq 2, \end{aligned} \tag{15}$$

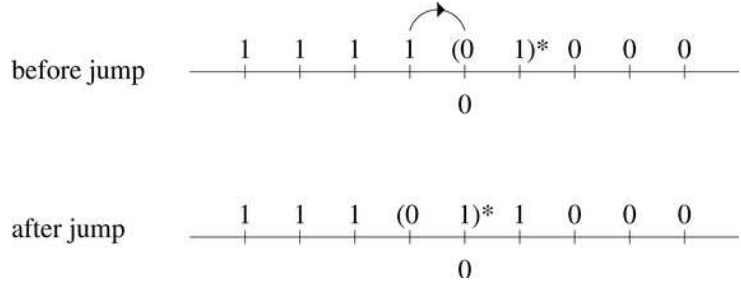

Fig. 3. *Pair representation of the second class particle.*

be the positions of the particles of $\eta_0^{01}$ labeled from right to left and the positions of the holes, labeled from left to right. We construct $P_i(t)$ and $H_i(t)$, the position of the $i$th particle, respectively $i$th hole, at time $t$ as a function of the random variables $G^{01}(z) := G(z) - w(1,1)$. The rule is:

(16)   at time $G^{01}(i,j)$ the $j$th particle and the $i$th hole interchange positions.

The initial ordered labels of the holes and particles make that after the $(j-1)$st particle has interchanged positions with the $i$th hole and the $j$th particle has interchanged positions with the $(i-1)$st hole, the $j$th particle must wait an exponential time of parameter 1 to interchange positions with the $i$th hole. This is the particle–hole interpretation of the recurrence relation (2).

Rule (16) is well defined in this case because only a finite number of exponential random variables is involved in the definition of each next move. Indeed, the variables $G^{01}(z)$ are well ordered, inducing a (random) order on the sites of $(\mathbb{Z}^+)^2$, say $z_1, z_2, \ldots$ with $G^{01}(z_k) < G^{01}(z_{k+1})$. In particular $z_1 \in \{(1,2),(2,1)\}$, for example. Starting with the minimum between $G(1,2)$ and $G(2,1)$, say $G(1,2) < G(2,1)$, then $z_1 = (1,2)$ and at time $G^{01}(z_2)$ the second particle and the first hole interchange positions (see Figure 4 ignoring the parentheses and the stars). Inductively, if $z_n = (i,j)$, then at time $G^{01}(z_n)$, the $j$th particle and the $i$th hole interchange positions. Call $P_i(G^{01}(z_n))$ and $H_i(G^{01}(z_n))$ the positions at time $G^{01}(z_n)$ of the $i$th particle and hole, respectively. For $i \geq 1$ define

$$
(17) \quad (P_i(t), H_i(t)) = (P_i(G^{01}(z_n)), H_i(G^{01}(z_n)))
$$
$$
\text{if } t \in [G^{01}(z_n), G^{01}(z_{n+1})).
$$

The resulting process $((P_i(t), H_i(t)), i \geq 1, t \geq 0)$ is the exclusion process in the sense that, if one disregards the labels, the process $(\zeta_t^{01}, t \geq 0)$ defined by

$$
(18) \quad \zeta_t^{01}(P_i(t)) = 1, \qquad \zeta_t^{01}(H_i(t)) = 0, \qquad i \geq 1,
$$

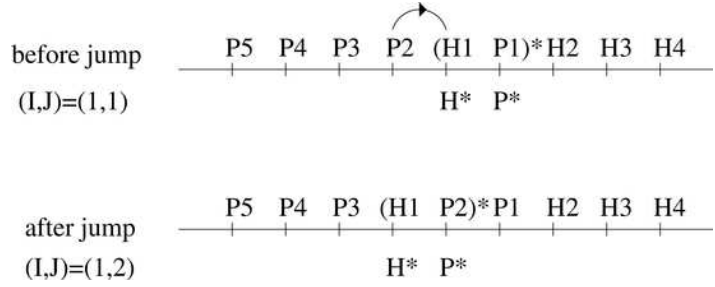

FIG. 4.  *Labels of particles, holes and \*pair. The particle configuration and the jump are the same as in Figure* 3.

has the same law as the process $(\eta_t^{01}, t \geq 0) = \Phi(\eta^{01}, \mathcal{N})$, defined with the Poisson processes. We call $\Upsilon(\eta^{01}, \mathcal{W}) = (\zeta_t^{01}, t \geq 0)$ the deterministic function that constructs $\zeta_t^{01}$ using $\mathcal{W}$.

*The second class particle in the competition model.* In the previous paragraph we have constructed a simple exclusion process starting with a particle at site 1 and a hole at site 0. In this construction we keep track of the position of each particle and hole. We now want to track the \*pair, the \*hole and \*particle initially at sites 0 and 1, respectively, whose evolution is described after Theorem 2. The labels of the \*particle and \*hole change with time. At time 0 the \*particle has label 1 and so does the \*hole: $P^*(0) = P_1(0)$ and $H^*(0) = H_1(0)$ and hence the labels of the \*pair are represented by the point $\varphi_0 = (1, 1)$, the initial value of the competition interface. Suppose in the next step, say $G(1, 2) < G(2, 1)$, the second particle jumps over the \*hole before the \*particle jumps over the second hole (see Figure 4). In this case, the labels of the \*pair at time $G(1, 2)$ are $(1, 2)$, which is exactly the argument that minimizes $\{G(2, 1), G(1, 2)\}$, so that, after the first jump of the \*pair, its labels are given by $\varphi_1$ [recall (8)]. By recurrence, $\varphi_n$ gives exactly the labels of the \*pair after its $n$th jump. More precisely, let $\tau_0 := 0$ and define

$$\tau_n := G^{01}(\varphi_n) \tag{19}$$

where $(\varphi_n, n \geq 0)$ is the competition interface defined in (7). The labels of the \*pair are given by the coordinates of the competition interface:

$$(H^*(\tau_n), P^*(\tau_n)) = (H_{i_n}(\tau_n), P_{j_n}(\tau_n)) \tag{20}$$

where $i_n$ and $j_n$ are the coordinates of $\varphi_n$: $(i_n, j_n) := \varphi_n$. Define the process $(\psi_t, t \geq 0) = (I(t), J(t), t \geq 0) \in (\mathbb{Z}^+)^2$ by

$$\psi_t := \varphi_n \qquad \text{if } t \in [\tau_n, \tau_{n+1}). \tag{21}$$

By definition (3) of $\gamma_t$, it is clear that $\psi_t$ belongs to both the growth interface and the competition interface (see Figure 1):

$$(22) \quad \psi_t \in \varphi \cap \gamma_{t+w(1,1)} \quad \text{and} \quad (H^*(t), P^*(t)) = (H_{I(t)}(t), P_{J(t)}(t)).$$

On the other hand, when the *pair jumps to the right the *hole increments its label by one, and when the *pair jumps to the left, the *particle increments its label by one. Hence,

$$(23) \qquad (H^*(t), P^*(t)) = (I(t) - J(t), I(t) - J(t) + 1).$$

Combining (23) with (14), we get the following result.

PROPOSITION 3. *The processes $((\bar{\eta}_t, I(t) - J(t)), \ t \geq 0)$ and $((\eta_t^1, X(t)), t \geq 0)$ are identically distributed.*

We construct simultaneously both processes in such a way that they are identical almost surely. See Lemma 6.

Using the technology of geodesics and the ergodicity of the last-passage percolation model we prove in Section 3 the following proposition (this is Theorem 1 without identifying the limit).

PROPOSITION 4.

$$(24) \qquad \lim_{t \to \infty} \frac{\varphi_n}{|\varphi_n|} = e^{i\theta}, \qquad \mathbb{P}\text{-}a.s.$$

where $\theta = \theta(\mathcal{W})$ is a random angle in $[0, 90°]$.

Propositions 3 and 4 and (22) are the keys to characterize the long time behavior of $(\psi_t, t \geq 0)$ as a line with a random angle and identify the distribution of the limiting angle:

PROPOSITION 5. *The following limits hold $\mathbb{P}$-a.s.:*

$$(25) \qquad \lim_{t \to \infty} \frac{\psi_t}{|\psi_t|} = e^{i\theta},$$

$$(26) \qquad \lim_{t \to \infty} \frac{\psi_t}{t} = e^{i\theta}/\mu(e^{i\theta}),$$

$$(27) \qquad \lim_{t \to \infty} \frac{I(t) - J(t)}{t} = f(\theta),$$

*where $\theta = \theta(\mathcal{W})$ is the random angle in $[0, 90°]$ given by Proposition 4,*

$$(28) \qquad f(\theta) := \frac{\sqrt{\cos \theta} - \sqrt{\sin \theta}}{\sqrt{\cos \theta} + \sqrt{\sin \theta}}$$

*and $f(\theta)$ is distributed uniformly in $[-1, 1]$:*

$$(29) \qquad \mathbb{P}(f(\theta) \leq u) = \tfrac{1}{2}(u + 1).$$

PROOF. Since $\psi_t \in \gamma_{t+w(1,1)}$ and by (6) $\inf\{|z| : z \in \gamma_t\}$ is of the order of $t$ [indeed, this infimum divided by $t$ converges to $1/\sqrt{8}$, the distance between the origin and the curve $\{\mu(u,v) = 1\}$], $|\psi_t| \to \infty$ as $t \to \infty$ and (25) follows from (24).

The limit (26) follows from (25), (22), the shape theorem (6) and (5). Indeed, the shape theorem (6) and the limit (25) imply that $\psi_t/t$ converges $\mathbb{P}$-almost surely to $g(\theta)e^{i\theta}$, where $g(\theta)$ is the distance from the origin to the intersection of the limiting curve $\mathcal{M} = \{(u,v) \in (\mathbb{R}^+)^2 : \mu(u,v) = 1\}$ with the line $\{(u,v) \in (\mathbb{R}^+)^2 : \tan\theta = u/v\}$ (the line with inclination $\theta$). Hence by the definition (5) of $\mu$, $\sqrt{g(\theta)\cos\theta} + \sqrt{g(\theta)\sin\theta} = 1$, from where (26) is derived.

The limit in (27) is an immediate consequence of (25) and (26). It is a uniform random variable as consequence of Proposition 3—that identifies the difference between the coordinates of the interface with the second class particle—and Ferrari and Kipnis [6], who proved that the asymptotic law of the second class particle is uniform in $[-1, 1]$. □

We finish this section by proving Theorems 1 and 2.

PROOF OF THEOREM 1. The $\mathbb{P}$-a.s. convergence follows from Proposition 4. Since by (29) $f(\theta)$ is uniformly distributed in $[-1, 1]$ and $f(\alpha)$ is decreasing in $\alpha$,

$$(30) \quad \mathbb{P}(\theta \leq \alpha) = \mathbb{P}(f(\theta) \geq f(\alpha)) = \frac{1}{2}(1 - f(\alpha)) = \frac{\sqrt{\sin\alpha}}{\sqrt{\sin\alpha} + \sqrt{\cos\alpha}}. \quad \Box$$

The proof of Theorem 2 requires the following lemma.

LEMMA 6. *There exists a map $R : \mathcal{N} \mapsto \mathcal{W}$ such that if the trajectory of the second class particle $(X(t), t \geq 0)$ as a function of $\mathcal{N}$ is well defined, then it is identical to the trajectory of $(I(t) - J(t), t \geq 0)$ as a function of $R(\mathcal{N})$. Furthermore, if $\mathcal{N}$ has law $\mathbb{Q}$, then $R(\mathcal{N})$ has law $\mathbb{P}$.*

PROOF. Let $\mathcal{N}$ be a family of Poisson processes. Let $((\eta_t^1, X(t)) : t \geq 0)$ be the exclusion process starting with the configuration full of particles to the left of the origin, empty to the right of the origin and with one second class particle in the origin constructed using $\mathcal{N}$.

Let $N$ be a Poisson process independent of $\mathcal{N}$. Let $\tau_n(\mathcal{N})$ be the times of jumps of the second class particle $X(t)$ with $\tau_0 = 0$. Then define $\mathcal{N}' = (N'_x(t) : t \geq 0)$ as a function of $\mathcal{N}$ and $N$ as follows:

$$(31) \quad N'_x[\tau_n, \tau_{n+1}) := \begin{cases} N_x(\tau_n, \tau_{n+1}], & \text{if } x < X(\tau_n), \\ N(\tau_n, \tau_{n+1}], & \text{if } x = X(\tau_n), \\ N_{x-1}(\tau_n, \tau_{n+1}], & \text{if } x > X(\tau_n). \end{cases}$$

Here $N_x(s,t]$ is the Poisson process $N_x$ in the interval $(s,t]$ (as a counting measure), and analogously for $N$. By the strong Markov property, $\mathcal{N}'$ has the same law as $\mathcal{N}$.

Let $\eta^{01}$ be the configuration defined in (13). Label its particles as in (15). Let the *pair be the *hole and the *particle initially at sites 0 and 1, respectively. Realize the process $\eta_t^{01}$ as a function of $\mathcal{N}'$. For this evolution track the position of the labeled particles $P_i(t)$ and holes $H_i(t)$ and the *pair $(H^*(t), P^*(t))$ as a function of the particle jumps as described after display (13). In this way we construct the processes $(\eta_t^{01}; P_i(t), H_i(t), i \geq 1; H^*(t), P^*(t); t \geq 0)$ as a function of $\mathcal{N}'$. Call $(I(t), J(t))$ the labels of the *hole–*particle at time $t$, so that $(H^*(t), P^*(t)) = (H_{I(t)}(t), P_{J(t)}(t))$; of course these are also function of $\mathcal{N}'$.

Then, for all $t$:

(32) $$X(t)(\mathcal{N}) = H_{I(t)}(\mathcal{N}') = I(t)(\mathcal{N}') - J(t)(\mathcal{N}'),$$

that is, the second class particle in the system governed by $\mathcal{N}$ is in the same place as the *hole in the system governed by $\mathcal{N}'$. Collapsing the *hole–*particle in the system governed by $\mathcal{N}'$, one obtains the particle configuration of the system governed by $\mathcal{N}$:

(33) $$\eta_t^{01}(\mathcal{N}')(x) := \begin{cases} \eta_t^1(\mathcal{N})(x), & \text{if } x < X(t)(\mathcal{N}), \\ 0, & \text{if } x = X(t)(\mathcal{N}), \\ 1, & \text{if } x = X(t)(\mathcal{N}) + 1, \\ \eta_t^1(\mathcal{N})(x-1), & \text{if } x > X(t)(\mathcal{N}) + 1. \end{cases}$$

Define $G'(1,1) = 0$ and for $i,j \geq 1$ let $G'(i,j) = G'(\mathcal{N}')(i,j)$ be the time the $j$th $\eta^{01}$ particle jumps over the $i$th hole. Define $w'(\mathcal{N}')(1,1)$ as an exponential random variable independent of $\mathcal{N}'$ and $w' = w'(\mathcal{N}')$ by $w'(i,j) = G'(i,j) - \max\{G'(i-1,j), G'(i,j-1)\}$. Since $w'(i,j)$ is the time the $i$th particle waits to jump over the $j$th hole when they are neighbors, $w'(i,j)$ are independent and identically distributed exponential of rate 1 (again strong Markov property). Hence $R(\mathcal{N}) := \mathcal{W}' = (w'(\mathcal{N}')(i,j) : \{(i,j) \in \mathbb{N}^2\})$ has the same law as $\mathcal{W}$.

It is immediate to check that

(34) $$(I(t), J(t))(\mathcal{W}') = (I(t), J(t))(\mathcal{N}') \quad \text{for all } t \geq 0.$$

That is, the *pair evolution described after (22) using the exponential times $\mathcal{W}'$ is exactly the same as the *pair evolution constructed as a function of the Poisson processes $\mathcal{N}'$. Notice that the auxiliary Poisson process $N$ used in the definition (31) of $\mathcal{N}'$ as a function of $\mathcal{N}$ plays no role in the *pair evolution. This is also true for $w'(1,1)$. □

PROOF OF THEOREM 2. The convergence $\mathbb{P}$-a.s. is established in Proposition 5. The convergence $\mathbb{Q}$-a.s. is a consequence of Lemma 6. □

**3. Geodesics.** In this section we prove Proposition 4. We introduce the notion of geodesics in last-passage percolation and explore its connection with the competition interface. Let $\pi = (z_k;\ k = 1, \ldots, n)$ be an up/right path from $z$ to $z'$. We say that $\pi$ is a *geodesic* from $z$ to $z'$ if

$$G(z, z') = \sum_{z'' \in \pi} w(z''). \tag{35}$$

Of course this is not a "geodesic" in the sense that it is the shortest way between two points. Indeed our geodesic is the longest oriented path between two points. For all $z, z' \in \mathbb{Z}$ there exists $\mathbb{P}$-a.s. a unique geodesic from $z$ to $z'$ which is denoted by $\pi(z, z')$. If $u = (u_1, u_2)$ and $v = (v_1, v_2)$ belongs to $\mathbb{R}^2$ and $u_k \leq v_k$ for $k = 1, 2$, then we define $G(u, v) = G(z_u, z_v)$ where $u \in Q(z_u)$ and $v \in Q(z_v)$, where we recall $Q(z)$ is the unit square with north-east point $z$. Analogously, we define $\pi(u, v) = \pi(z_u, z_v)$. Let $\pi_z = (z_k;\ k = 1, \ldots)$ be an *up/right semi-infinite path* starting at $z = z_1$. For each $\alpha \in [0, 90°]$ we say that $\pi_z$ has direction $\alpha$ if

$$\lim_{k \to \infty} \frac{z_k}{|z_k|} = e^{i\alpha}.$$

We say that $\pi_z$ is a *uni-geodesic* if for all $i < j$ the geodesic from $z_i$ to $z_j$ is exactly $(z_i, \ldots, z_j)$. For each $\alpha \in [0, 90°]$ we say that $\pi_z$ is an $\alpha$-*geodesic* if it is a uni-geodesic and has direction $\alpha$. Proposition 4 is a consequence of the following propositions concerning geodesics. The proofs follow Newman [15], Licea and Newman [10] and Howard and Newman [9] who proved analogous results for two-dimensions first-passage percolation models. Martin [13] has independently proved these results for the model under consideration.

Let $z \in \mathbb{Z}^2$ and $\mathbb{N}_z^2 = z + \mathbb{N}^2$ ($\mathbb{N} = \{0, 1, 2, \ldots\}$). Define $R(z) = \bigcup_{z' \in \mathbb{N}_z^2} \pi(z, z')$. Since $\mathbb{P}$-a.s. finite geodesics do exist and are unique, $R(z)$ can be seen as a tree spanning all $\mathbb{N}_z^2$. The set of vertices of the tree is $\mathbb{N}_z^2$ and the set of edges is $\{(z'', z') : z'' - z' = 1 \text{ and } z'' \in \pi(z, z')\}$.

PROPOSITION 7. *For $z \in \mathbb{Z}^2$ let $\Omega_1(z)$ be the event "every uni-geodesic $\pi_z \subseteq R(z)$ is an $\alpha$-geodesic for some $\alpha = \alpha(\pi_z) \in [0, 90°]$ and there exists at least one $\alpha$-geodesic for each $\alpha \in [0, 90°]$." Then $\mathbb{P}(\Omega_1(z)) = 1$.*

PROPOSITION 8. *For $z \in \mathbb{Z}^2$ and $\alpha \in [0, 90°]$ let $\Omega_2(z, \alpha)$ be the event "there exists at most one $\alpha$-geodesic in $R(z)$" and let $\ell$ be the Lebesgue measure in $[0, 90°]$. Then there exists a set $D \subseteq [0, 90°]$ of full Lebesgue measure such that for all $\alpha \in D$, $\mathbb{P}(\Omega_2(z, \alpha)) = 1$.*

We recall that $D$ does not depend on the realization of the exponential times $\mathcal{W}$. Indeed, a stronger version of Proposition 8 holds: for *every* $\alpha \in (0, 90°)$ there is only one $\alpha$-geodesic in $R(z)$ with probability 1 [13]. On the other hand, with probability 1 there are directions with more than one geodesic: $\mathbb{P}(\bigcap_{\alpha \in (0, 90°)} \Omega_2(z, \alpha)) = 0$.

For $\alpha \in D$ let $\pi_z(\alpha)$ be the unique $\alpha$-geodesic starting at $z$. This is $\mathbb{P}$-a.s. well defined by Propositions 7 and 8.

PROPOSITION 9. *For $\alpha \in D$ let $\Omega_3(\alpha)$ be the event "for all $z, z' \in \mathbb{Z}^2$, there exists a random point $c_\alpha = c(\alpha, z, z') \in \mathbb{Z}^2$ such that $\pi_z(\alpha) = \pi(z, c_\alpha) \cup \pi_{c_\alpha}(\alpha)$ and $\pi_{z'}(\alpha) = \pi(z', c_\alpha) \cup \pi_{c_\alpha}(\alpha)$." Then $\mathbb{P}(\Omega_3(\alpha)) = 1$.*

As a consequence of the above propositions we get that for all $\alpha \in D$, $\mathbb{P}$-a.s. for all $z, z' \in \mathbb{Z}^2$, $z \neq z'$, there exists a random $c_\alpha = c(\alpha, z, z')$ and $r_0 > 0$ such that for all $r > r_0$

$$(36) \qquad G(z, re^{i\alpha}) - G(z', re^{i\alpha}) = G(z, c_\alpha) - G(z', c_\alpha) \neq 0.$$

Indeed, from Propositions 7 and 8, if we fix $\alpha \in D$, then $\mathbb{P}$-a.s. for all $z \in \mathbb{Z}^2$ $\lim_{r \to \infty} \pi(z, re^{i\alpha}) = \pi_z(\alpha)$. This means that for all $\bar{z} \in \pi_z(\alpha)$ there exists $r_0 > 0$ such that for all $r > r_0$, $\pi(z, \bar{z}) \subseteq \pi(z, re^{i\alpha})$. This together with Proposition 9 implies that for all $z, z' \in \mathbb{Z}^2$ there exists $c_\alpha \in \mathbb{Z}^2$ and $r_0 > 0$ such that for all $r > r_0$, $\pi(z, re^{i\alpha}) = \pi(z, re^{i\alpha}) \cup \pi(c_\alpha, re^{i\alpha})$ and $\pi(z', re^{i\alpha}) = \pi(z', re^{i\alpha}) \cup \pi(c_\alpha, re^{i\alpha})$, which yields (36).

PROOF OF PROPOSITION 4. Let $\overline{\mathbf{G}}_\infty^{21} := \bigcup_{z \in \mathbf{G}_\infty^{21}} Q(z)$, $\overline{\mathbf{G}}_\infty^{12} := \bigcup_{z \in \mathbf{G}_\infty^{12}} Q(z)$. For each $\alpha \in [0, 90°]$ and $r > 0$ let $l_r^\alpha = \{se^{i\alpha};\ s > r\}$. Define the random decomposition of $[0, 90°]$ by

$$I_{21} = \{\alpha \in [0, 90°];\ \exists r_0 \text{ so that } l_{r_0}^\alpha \subseteq \overline{\mathbf{G}}_\infty^{21}\},$$

$$I_{12} = \{\alpha \in [0, 90°];\ \exists r_0 \text{ so that } l_{r_0}^\alpha \subseteq \overline{\mathbf{G}}_\infty^{12}\},$$

and $I = (I_{21} \cup I_{12})^c$. Notice that $0 \in I_{21}$, $90° \in I_{12}$, and since $\overline{\mathbf{G}}_\infty^{21}, \overline{\mathbf{G}}_\infty^{12}$ are connected regions of $\{(x, y);\ x > 0, y > 0\}$, then $I_{21}$ and $I_{12}$ are intervals in $[0, 90°]$. This implies that $I$ is also an interval in $[0, 90°]$. Thus if we denote $\varphi_n = |\varphi_n|e^{i\theta_n}$, then

$$(37) \qquad \left(\liminf_n \theta_n, \limsup_n \theta_n\right) \subseteq I.$$

Let $D_0$ be an enumerable subset of $D$ that is dense in $(0, 90°)$ (recall that $D$ has full Lebesgue measure). By (36), $\mathbb{P}$-a.s., for all $\alpha \in D_0$,

$$(38) \quad \begin{aligned} \lim_{r \to \infty} (G((2,1), re^{i\alpha}) &- G((1,2), re^{i\alpha})) \\ &= G((2,1), c_\alpha) - G((1,2), c_\alpha) \neq 0. \end{aligned}$$

Notice also that if $\alpha \in I$, then

$$
\begin{aligned}
(39) \quad & \liminf_{r \to \infty}(G((2,1), re^{i\alpha}) - G((1,2), re^{i\alpha})) \\
& \leq 0 \leq \limsup_{r \to \infty}(G((2,1), re^{i\alpha}) - G((1,2), re^{i\alpha})),
\end{aligned}
$$

because the line $l_0^\alpha$ alternates infinitely often its color and this implies (39). Thus (38) and (39) imply that

$$
(40) \qquad\qquad\qquad \mathbb{P}(I \cap D_0 = \varnothing) = 1.
$$

Now, (40) implies that $\mathbb{P}$-a.s. $I$ has empty topological interior and this together with (37) implies that $(\theta_n)_{n \in \mathbb{N}}$ converges. $\square$

The following lemma, proven in the end of this section, is the main ingredient to prove Proposition 7. It gives an upper bound for the fluctuations of the geodesics. Let $d(z, A)$ be the Euclidean distance between $z \in \mathbb{R}^2$ and the set $A \subset \mathbb{R}^2$.

LEMMA 10. *There exists $\varepsilon_0 > 0$ such that for all $\varepsilon \in (0, \varepsilon_0)$, there exist constants $C_1, C_2, C_3 > 0$ and $\delta > 0$ such that for all $z \in \mathbb{N}^2$ with $|z| > C_1$,*

$$
\mathbb{P}\bigg(\sup_{z' \in \pi(0,z)} d(z', [0,z]) \geq |z|^{3/4+\varepsilon}\bigg) \leq C_2 \exp(-C_3 |z|^\delta).
$$

PROOF OF PROPOSITION 7. By translation invariance we can assume $z = (0,0)$. For $\bar{z}, z' \in \mathbb{N}^2 \setminus \{(0,0)\}$, denote by $\text{ang}(\bar{z}, z')$ the angle in $[0, 90°]$ between $\bar{z}$ and $z'$ and let $C(\bar{z}, \varepsilon) = \{z';\ \text{ang}(\bar{z}, z') \leq \varepsilon\}$. Let $R$ be an infinite connected tree with vertices in $\mathbb{N}^2$ and nearest-neighbor oriented edges. Assume also that $(0,0)$ and $\bar{z} \in \mathbb{N}^2$ are vertices of $R$. We denote by $R^{\text{out}}[\bar{z}]$ the set of vertices $z'$ of $R$ such that the path in $R$ between $(0,0)$ and $z'$ touches $\bar{z}$. Let $h : \mathbb{R}^+ \to \mathbb{R}^+$. We say that $R$ is $h$-straight if for all but finitely many vertices $\bar{z}$ of $R$, $R^{\text{out}}[\bar{z}] \subseteq C(\bar{z}, h(|\bar{z}|))$. By Proposition 2.8 of [9], if $R$ is $h$-straight with $h$ satisfying $\lim_{L \to \infty} h(L) = 0$, then every semi-infinite path in $R$ starting from $(0,0)$ has a direction $\alpha \in [0, 90°]$ and for every $\alpha \in [0, 90°]$ there exists at least one semi-infinite path in $R$ starting from $(0,0)$ and with direction $\alpha$. Let $\delta \in (0,1)$ and set $h_\delta(L) = L^{-\delta}$. By Lemma 2.7 of [9], to prove that for all $\delta \in (0, 1/4)$, $R((0,0))$ is $h_\delta$-straight it is sufficient to prove that for all sufficiently small $\varepsilon > 0$, the number of $z \in \mathbb{N}^2$ such that $\sup_{z' \in \pi(0,z)} d(z', [0,z]) \geq |z|^{3/4+\varepsilon}$ is $\mathbb{P}$-a.s. finite. Therefore, by Borel–Cantelli, Proposition 7 is a consequence of Lemma 10. $\square$

PROOF OF PROPOSITION 8. Again we can assume that $z = (0,0)$. Let $e = (z, z+(1,0))$ be an edge of the tree $R((0,0))$ such that $z + (1,0)$ has infinitely many descendants. We inductively define a uni-geodesic $\pi_e$ in $R((0,0))$ as follows. Put $z_0 = z$ and $z_1 = z + (1,0)$. For each $n \geq 1$, if $z_n$ has exactly one child, say $z'_n$, with infinitely many descendants, then put $z_{n+1} = z'_n$; otherwise put $z_{n+1} = z_n + (0,1)$. If there are two distinct $\alpha$-geodesics starting from $(0,0)$, say $\pi_1$ and $\pi_2$, then they have to bifurcate at some $z \in R((0,0))$ going, respectively, to $z + (1,0)$ and $z + (0,1)$ in their next steps. In this case, $\pi_e$ with $e = (z, z+(1,0))$ is caught between $\pi_1$ and $\pi_2$. Hence $\pi_e$ is an $\alpha$-geodesic because we are in two dimensions. Therefore, $\Omega_2((0,0), \alpha)$ must occur unless the event $B(e, \alpha) := [\pi_e \text{ is an } \alpha\text{-geodesic}]$ occurs for some $e = (z, z+(1,0))$. Thus

$$(41) \qquad 1 \geq \mathbb{P}(\Omega_2((0,0), \alpha)) \geq 1 - \sum_{e=(z,z+(1,0))} \mathbb{P}(B(e, \alpha)).$$

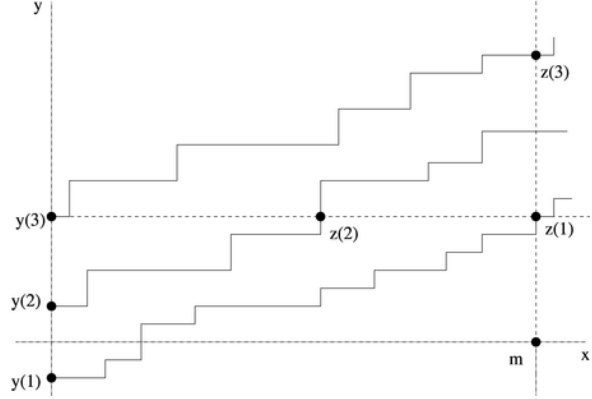

Fig. 5. *Local modification.*

For each $e = (z, z+(1,0))$, $\pi_e$ cannot be an $\alpha$-geodesic for more than one $\alpha$, and so $\int \mathbb{1}_{B(e,\alpha)} \ell(d\alpha) = 0$ for each realization of the exponential times. By Fubini,

$$(42) \qquad \int \mathbb{P}(B(e,\alpha))\ell(d\alpha) = \int \left[\int \mathbb{1}_{B(e,\alpha)} \ell(d\alpha)\right] d\mathbb{P} = 0.$$

Integrating (41) with respect to $\ell(d\alpha)$ and using (42) completes the proof of Lemma 8. □

PROOF OF PROPOSITION 9. By Proposition 8, for fixed $\alpha \in D$ $\mathbb{P}$-a.s., if $\pi_z(\alpha)$ and $\pi_{z'}(\alpha)$ are not site disjoint, then they must coalesce. Therefore we must show that there is zero probability that there exist disjoint $\alpha$-geodesics. Let $S(\alpha) = \bigcup_{z \in \mathbb{Z}^2} \pi_z(\alpha)$ be the set of $\alpha$-geodesics emanating from $z \in \mathbb{Z}^2$. Then $S(\alpha)$ is a forest composed by a random number $N(\alpha) \in \{1, \ldots, +\infty\}$ of connected trees. The event "there are not disjoint $\alpha$-geodesics" is equivalent to the event "$N(\alpha) = 1$" and with this formulation we can apply the Burton and Keane [2] argument. This argument is based on a local modification idea that is formalized as follows. Let $y_1 < \cdots < y_k$ be points in $\mathbb{Z}$, and let $A(y_1, \ldots, y_k)$ be the event "$\pi_{(0,y_1)}(\alpha), \ldots, \pi_{(0,y_k)}(\alpha)$ are disjoint and every site touched by $\pi_{(0,y_j)}$ after its initial site at $(0, y_j)$ has strictly positive first coordinate" (see Figure 5). We claim that

(43)　　if $\mathbb{P}(N(\alpha) \geq 2) > 0 \implies \exists y_1, y_2, y_3; \mathbb{P}(A(y_1, y_2, y_3)) > 0.$

Indeed, if the right-hand side of (43) holds, then there exist $y_1, y_2$ such that the probability of $A(y_1, y_2)$ is positive. For $i = 1, 2$ and $l \in \mathbb{Z}$ let $y_i^l = y_i + l(y_1 - y_2)$. By translation invariance, the probability of $A(y_1^l, y_2^l)$ does not depend on $l$. This together with Fatou implies that for some $l_1 < l_2$ the probability of $A(y_1^{l_1}, y_2^{l_2}) \cap A(y_1^{l_2}, y_2^{l_2})$ is also positive, although $\pi_{y_1^{l_1}}$ cannot intersect $\pi_{y_1^{l_2}}$ or $\pi_{y_2^{l_2}}$ because otherwise it must intersect $\pi_{y_2^{l_1}}$ (by planarity). Thus $A(y_1^{l_1}, y_2^{l_2}) \cap A(y_1^{l_2}, y_2^{l_2}) \subseteq A(y_1^{l_1}, y_1^{l_2}, y_2^{l_2})$ which proves (43). Let $A_m^\delta(y_1, y_2, y_3) := A(y_1, y_2, y_3) \cap B_\delta \cap C_m$, where $B_\delta$ is the event "$w(z) > \delta$ for all $z \in [(-1, y_1), (-1, y_3)] \cup [(-1, y_3), (m, y_3)]$" and $C_m$ is the event "$\pi_{(0,y_1)}(\alpha)$ intersects the point $(m, y_3)$." Since $\alpha \in (0, 90°)$, (43) implies that, if the right-hand side of (43) holds, then for some $m > 0$ and $\delta > 0$,

(44)　　　　　　　　$\mathbb{P}(A_m^\delta(y_1, y_2, y_3)) > 0.$

Consider the event $A = A_m^\delta(y_1, y_2, y_3)$ appearing in (44). Let $z_2$ be the first intersection between $\pi_{(0,y_2)}(\alpha)$ and $[(0, y_3), (m, y_3)]$, and let $\pi_1$ be the piece of $\pi_{(0,y_1)}(\alpha)$ between the points $(0, y_1)$ and $z_1 := (m, y_3)$ (see Figure 5). Consider the bounded region $\Lambda$ (not including the boundary) limited by $[(-1, y_1), (-1, y_3)]$, $[(-1, y_3), (m, y_3)]$ and $\pi_1$. Thus we define a mapping $\Phi$ on subsets $B$ of $A$ by first letting $W(\omega) = \{z \in \Lambda; w(z) > \delta\}$ and then setting

$$\Phi(B) = \bigcup_{\omega \in B} \left[ \prod_{z \notin W(\omega)} \{w(z)\} \times \prod_{z \in W(\omega)} (0, \delta) \right].$$

Heuristically, $\Phi$ alters each $\omega \in B$ into an $\omega' \in \Phi(B)$ by changing each $w(z) > \delta$ with $z \in \Lambda$ to some value $w'(z) \in (0, \delta)$ [it may happen that $\Phi(B)$ is nonmeasurable]. Since the $w(z)$'s with $z \in \pi_{(0,y_1)}(\alpha)$, or $z \in \pi_{(0,y_3)}(\alpha)$, or $z \in \pi_{z_2}(\alpha)$, were unchanged while the others decreased or stayed as before, it follows that each one of the paths $\pi_{(0,y_1)}(\alpha)$, $\pi_{(0,y_3)}(\alpha)$ and $\pi_{z_2}(\alpha)$ continues to be an $\alpha$-geodesic for $\omega' \in \Phi(A)$. Similarly, $\omega'$ continues to belong to $B_\delta$. Although, since for each $z \in \Lambda$ we have $w'(z) < \delta$ and for each $z \in [(-1, y_1), (-1, y_3)] \cup [(-1, y_3), (m, y_3)]$ we have $w'(z) > \delta$, then for all $z \in [(-1, y_1), (-1, y_3)]$ and $z' \in [(-1, y_3), (m, y_3)]$ the geodesic $\pi(z, z')$ for $\omega'$ either will be the path which starts at $z$, goes vertically until it reaches $(-1, y_3)$ and then goes horizontally until it reaches $z'$, or $\pi(z, z') \cap \pi_1 \neq \varnothing$. Let $z_3 := (m, y')$ be the point where the geodesic from $(0, y_3)$ crosses the vertical line $\{m\} \times \mathbb{Z}$. Therefore, any $\alpha$-geodesic for $\omega'$, starting at $u = (x, y) \notin [0, m] \times [y_1, y']$ with $x < m$, cannot touch the middle path $\pi_{z_2}(\alpha)$ without first intersecting $\pi_{(0,y_1)}$ or $\pi_{(0,y_3)}$ and this leads to a contradiction because in such a case, by Proposition 8, they must coalesce. Thus, $\Phi(A) \subseteq F$ where $F$ denotes the event that some tree in $S(\alpha)$ touches the rectangle $[0, m] \times [y_1, y']$ but no other site in the half-plane $\{(x, y);\ x < m\}$. Since to each site we attached an exponential random variable and $\mathbb{P}(A) > 0$ [inequality (44)], by Lemma 3.1 of [10], there exists a measurable set $\bar{A} \subseteq \Phi(A)$ such that $\mathbb{P}(\bar{A}) > 0$, which implies that $\mathbb{P}(F) > 0$. Now consider a rectangular array of nonintersecting translates $\Theta_z$ of the rectangle $\Theta_0 = [0, m] \times [y_2, y_3]$ indexed by $z \in \mathbb{Z}$, and consider the corresponding translated events $F_z$ of $F_0 = F$. Notice that, if $F_z$ and $F_{z'}$ both occur, with $z \neq z'$, then the corresponding trees in $S(\alpha)$ must be disjoint. Let $n_L$ be the number of $\Theta_z$'s in $[0, L]^2$ and let $N_L$ be the number of the corresponding $F_z$'s which occur. By translation invariance $\mathbb{E} N_L = n_L \mathbb{P}(F)$. The number of disjoint trees in $S(\alpha)$ which touch $[0, L]^2$ cannot exceed the number of boundary sites in $[0, L]^2$ and this together with $\mathbb{P}(F) > 0$ yields a contradiction for large $L$ because $n_L$ is of order $L^2$. Therefore, we have proved that $\mathbb{P}(F) = 0$ and this together with (43) and (44) implies that $\mathbb{P}(N(\alpha) \geq 2) = 0$. This together with Proposition 7 implies that $\mathbb{P}(N(\alpha) = 1) = 1$ which completes the proof of Proposition 9. □

The proof of Lemma 10 is based on the following lemma that provides an upper bound for moderate deviations of $G(0, z)$ from its asymptotic value $\mu(z)$.

LEMMA 11. *There exists $\varepsilon_0 > 0$ such that for all $\varepsilon \in (0, \varepsilon_0)$, there exist constants $C_4, C_5, C_6 > 0$ such that for all $z \in \mathbb{N}^2$ with $|z| > C_4$, for any $r \in [|z|^{1/2+\varepsilon}, |z|^{3/2-\varepsilon}]$*

$$\mathbb{P}(|G(0, z) - \mu(z)| \geq r) \leq C_5 \exp(-C_6 r / |z|^{1/2}).$$

We prove this lemma after the proof of Lemma 10.

PROOF OF LEMMA 10. Let $C_\varepsilon(z) = \{z';\ d(z', [0, z]) \geq |z|^{3/4+\varepsilon}\}$ and let $\Delta(z, z') = \mu(z) - \mu(z - z') - \mu(z')$. If $\sup_{z' \in \pi(0, z)} d(z', [0, z]) \geq |z|^{3/4+\varepsilon}$, then there exists $z' \in \partial C_\varepsilon(z)$ such that $G(0, z) = G(0, z') + G(z', z)$. By summing $\Delta(z, z')$ in both sides of the last equality and using the translation invariance of the model, we obtain

$$\mathbb{P}\left(\sup_{z' \in \pi(0, z)} d(z', [0, z]) \geq |z|^{3/4+\varepsilon}\right)$$

$$
\begin{align}
(45) \quad &\leq \mathbb{P}(|G(0,z) - \mu(z)| \geq |\Delta(z,z')|/3) \\
&+ \mathbb{P}(|G(0,z') - \mu(z')| \geq |\Delta(z,z')|/3) \\
&+ \mathbb{P}(|G(0,z-z') - \mu(z-z')| \geq |\Delta(z,z')|/3).
\end{align}
$$

If $z = (z_1, z_2), z' = (z'_1, z'_2)$, then

$$\Delta(z, z') = 2(\sqrt{z_1 z_2} - \sqrt{z'_1 z'_2} - \sqrt{(z_1 - z'_1)(z_2 - z'_2)}).$$

This implies that (see Lemma 2.1 of [23]) there exist constants $A_1, A_2, A_3 > 0$ such that for all $z \in \mathbb{N}^2$ with $|z| \geq A_1$, for all $z' \in \partial C_\varepsilon(z)$

$$(46) \qquad A_2 |z|^{1/2 + 2\varepsilon} \leq \Delta(z, z') \leq A_3 |z|^{3/4 + \varepsilon}.$$

Notice that (46) implies that there exists $M > 0$, such that for all $z \in \mathbb{N}^2$ with $|z| > M$, $\Delta(z, z') \in [\bar{z}^{1/2 + 2\varepsilon}, \bar{z}^{3/4 + \varepsilon}]$, where $\bar{z} = z$, or $\bar{z} = z'$, or $\bar{z} = z - z'$. By choosing $\varepsilon > 0$ small enough and using (45), together with Lemma 11 we complete the proof of Lemma 10. $\square$

PROOF OF LEMMA 11. Since for all $z = (z_1, z_2) \in \mathbb{N}^2$, for all $\pi \in \Pi(0, z)$ (the set of all *up-right* paths connecting 0 to $z$), $|\pi| = z_1 + z_2 + 1$, where $|\pi|$ is the number of sites in $\pi$, it is a consequence of Corollary 8.2.4 of [22] that there exist constants $A_1, A_2, A_3 > 0$ such that for all $z \in \mathbb{N}^2$, for all $x \in [0, A_1 |z|]$,

$$(47) \qquad \mathbb{P}(|G(0,z) - \mathbb{E}G(0,z)| > x\sqrt{|z|}) \leq A_2 \exp(-A_3 x).$$

To replace $\mathbb{E}G(0,z)$ by $\mu(0,z)$ in (47) we need to consider some technical details. First, we claim that (47) implies that there exist constants $A_4, A_5$ such that for all $z \in \mathbb{N}^2$ with $|z| \geq A_4$,

$$(48) \qquad \mathbb{E}G(0, 2z) \leq 2\mathbb{E}G(0, z) + A_5 \sqrt{|z|} \log(|z|).$$

Indeed, let $H_z = \{z'; |z'| = |z|\} \cap \mathbb{N}^2$ ($|z| = |z_1| + |z_2|$). By the definition of $G$ and by the translation invariance of the model,

$$(49) \qquad \mathbb{E}G(0, 2z) \leq 2\mathbb{E} \max_{z' \in H_z} G(0, z').$$

Now assume that $Y_i^z$ for $1 \leq i \leq n(z)$ are nonnegative random variables on a common probability space such that for some $a, M_0, C_0, C_1, C_2, C_3 \in (0, +\infty)$, for all $z \in \mathbb{N}^2$ with $|z| > M_0$,

$$(50) \qquad \mathbb{E}(Y_i^z) \leq |z|^a \quad \text{and} \quad n(z) \leq C_0 |z|,$$

and

$$(51) \qquad \mathbb{P}(|Y_i^z - \mathbb{E}(Y_i^z)| > x) \leq C_1 \exp(-C_2 x) \qquad \text{for } x \leq C_3 |z|.$$

Then, for some $M_1, C_4 > 0$ (see Lemma 4.3 of [9]), for all $z \in \mathbb{N}^2$ with $|z| > M_1$,

$$\mathbb{E} \max_{1 \leq i \leq n(z)} (\mathbb{E}(Y_i^z) - Y_i^z) \leq C_4 \log |z|. \tag{52}$$

Therefore, to conclude the proof of (48) we order the points $z' \in H(z)$ by $z_1, \ldots, z_n$, where $n$ depends on $z$ but $n \leq C_0 |z|$ for some constant $C_0$. Take $Y_i^z = G(0, z_i)/\sqrt{|z|}$ and note that the hypotheses (50), (51) are satisfied with $a = \varepsilon + 1/2$, $C_0$ as before, and $M_0, C_1, C_2, C_3$ given by (47). Thus, (52) together with (49) completes the proof of (48). Now, we claim that the superadditivity of $\mathbb{E}G(0, z)$ and (48) imply that for some constant $A_6 > 0$

$$\mathbb{E}G(0, z) - A_6 \sqrt{|z|} \log(|z|) \leq \mu(z) \leq \mathbb{E}G(0, z). \tag{53}$$

Indeed, the right-hand side of (53) is an immediate consequence of superadditivity. To prove the left-hand side, assume that $h: \mathbb{R}^+ \to \mathbb{R}$ and $g: \mathbb{R}^+ \to \mathbb{R}^+$ satisfy the following conditions: $\lim_{s \to \infty} h(s)/s\mu \in \mathbb{R}$, $\lim_{s \to \infty} g(s)/s = 0$, $h(2s) \geq 2h(s) - g(s)$ and $\phi = \limsup_{s \to \infty} g(2s)/g(s) < 2$. Then, for any $c > 1/(2 - \phi)$, $h(s) \leq \mu s + cg(s)$ for all large $s$ (see Lemma 4.2 of [9]). Therefore, if we fix a direction $\hat{z} = z/|z|$ and take $h(s) = -\mathbb{E}G(0, s\hat{z})$, $g(s) = A_5 \sqrt{s} \log(s)$, then this last claim together with (48) completes the proof of (53). Thus, (53) and (47) imply that for some constants $A_7, A_8, A_9, A_{10} > 0$, for all $x \in [A_7 \log |z|, A_8 |z|]$,

$$\mathbb{P}(|G(0, z) - \mu(z)| > 2x\sqrt{|z|}) \leq A_9 \exp(-A_{10} x). \tag{54}$$

Taking $r = 2x\sqrt{|z|}$ and adjusting the constants, (54) yields Lemma 11. □

**4. Final remarks.** We have shown a law of large numbers for the competition interface in last-passage percolation in the positive quadrant $(\mathbb{Z}^+)^2$. A crucial step in this proof was Proposition 9 which establishes that unigeodesics starting at different fixed points with the same direction must coalesce. The law of large numbers for the competition interface also holds for other random regions as a consequence of the law of large numbers for the second class particle of Mountford and Guiol [14] and Lemma 6. These regions are limited to the south-west by a random curve $\gamma = (\gamma_n, n \in \mathbb{Z}) \subset \mathbb{Z}^2$ defined by $\gamma_0 = (1, 1)$, $\gamma_1 = (1, 0)$, $\gamma_{-1} = (0, 1)$ and then $\gamma_n - \gamma_{n-1} = (\eta(n) - 1, -\eta(n))$, for $n \in \mathbb{Z} \setminus \{0, 1\}$ and $\eta$ distributed according to the product measure with densities $\lambda$ to the left of the origin and $\rho$ to the right of it. Since Lemma 6 can be extended to any region obtained as a transformation of the initial configuration of the simple exclusion process, the law of large numbers for the competition interface also holds in this case [7]. However, it would be nice to have an autonomous proof using geodesics. To extend the result to the regions considered by Mountford and Guiol one should be able to show that when the point is asymptotically beyond the corresponding characteristic the "point to semi-line" geodesic is realized in the limit by a random location in the semi-line. More precisely, let $L_\rho$ be a random semi-line starting at $(0,0)$ doing independent steps at right with probability $1-\rho$ and down with probability $\rho$. This interface corresponds to the right initial configuration for the simple exclusion process chosen with the product measure with density $\rho$. Let $z_n = (x_n, y_n)$ be a sequence of points in $\mathbb{N}^2$ such that $x_n, y_n \to \infty$ and $\frac{x_n}{y_n} \to (\frac{\rho}{1-\rho})^2 - \varepsilon$ for some $\varepsilon > 0$, as $n \to \infty$. Let $g_n$ be the location in $L_\rho$ that realizes the $z_n$ to $L_\rho$ geodesic. Then one needs to show that as $n \to \infty$, $g_n \to g$, a random location, almost surely. The inclination $(\frac{\rho}{1-\rho})^2$ corresponds to the asymptotic behavior of the second class particle under this initial measure: $\frac{X(t)}{t} \to (1-2\rho) = (1-\rho)^2 - \rho^2$, as $t \to \infty$.

An anonymous referee and Christoffe Bahadoran asked the authors about the resemblance between our Proposition 3 which identifies the second class particle and the competition interface determined by looking, in the last-passage picture, from which side of point $(1,1)$ the maximizing paths of different points emanate. The referee says: "This bears a curious resemblance to Proposition 4.1 in [21]: that result also identifies the position $X(t)$ of the second class particle by looking at which side of the initial position $X(0)$ come the maximizers in the variational formula of the process. One wonders whether these two representations are two sides of the same coin." We leave this investigation for future work.

**Acknowledgments.** We thank James Martin, Christoffe Bahadoran and Tom Mountford for enlightening discussions and a referee for a careful reading and useful comments about a previous version of this paper.

## REFERENCES

[1] BENASSI, A. and FOUQUE, J.-P. (1987). Hydrodynamical limit for the asymmetric simple exclusion process. *Ann. Probab.* **15** 546–560. MR885130
[2] BURTON, R. and KEANE, M. (1989). Density and uniqueness in percolation. *Comm. Math. Phys.* **121** 501–505. MR990777
[3] DERRIDA, B. and DICKMAN, R. (1991). On the interface between two growing Eden clusters. *J. Phys. A* **24** L191–193. Available at http://stacks.iop.org/JPhysA/24/L191.
[4] FERRARI, P. A. (1992). Shock fluctuations in asymmetric simple exclusion. *Probab. Theory Related Fields* **91** 81–101. MR1142763
[5] FERRARI, P. A. (1994). Shocks in one-dimensional processes with drift. In *Probability and Phase Transition* (G. Grimmett, ed.) **420** 35–48. Kluwer, Dordrecht. MR1283174
[6] FERRARI, P. A. and KIPNIS, C. (1995). Second class particles in the rarefaction front. *Ann. Inst. H. Poincaré* **31** 143–154. MR1340034
[7] FERRARI, P. A., MARTIN, J. B. and PIMENTEL, L. P. R. (2004). Roughening and inclination properties of competition interfaces. ArXiv:math. PR/0412198.

[8] HÄGGSTROM, O. and PEMANTLE, R. (1998). First-passage percolation and a model for competing spatial growth. *J. Appl. Probab.* **35** 683–692. [MR1659548](MR1659548)
[9] HOWARD, C. D. and NEWMAN, C. M. (2001). Geodesics and spanning trees for Euclidean first-passage percolation. *Ann. Probab.* **29** 577–623. [MR1849171](MR1849171)
[10] LICEA, C. and NEWMAN, C. M. (1996). Geodesics in two dimension first-passage percolation. *Ann. Probab.* **24** 399–410. [MR1387641](MR1387641)
[11] LIGGETT, T. M. (1985). *Interacting Particle Systems.* Springer, New York. [MR776231](MR776231)
[12] LIGGETT, T. M. (1999). *Stochastic Interacting Systems*: *Contact, Voter and Exclusion Processes.* Springer, New York. [MR1717346](MR1717346)
[13] MARTIN, J. (2004). Unpublished manuscript.
[14] MOUNTFORD, T. and GUIOL, H. (2005). The motion of a second class particle for the TASEP starting from a decreasing shok profile. *Ann. Appl. Probab.* **15** 1227–1259. [MR2134103](MR2134103)
[15] NEWMAN, C. M. (1995). A surface view of first-passage percolation. In *Proceedings of International Congress of Mathematicians 1994* (S. D. Chatterji, ed.) **2** 1017–1023. Birkhäuser, Basel. [MR1404001](MR1404001)
[16] PIMENTEL, L. P. R. (2004). Competing growth, interfaces and geodesics in first-passage percolation on Voronoi tilings. Ph.D. thesis, IMPA, Rio de Janeiro. Available at www.impa.br/preprint.
[17] REZAKHANLOU, F. (1991). Hydrodynamic limit for attractive particle systems on $Z^d$. *Comm. Math. Phys.* **140** 417–448. [MR1130693](MR1130693)
[18] REZAKHANLOU, F. (1995). Microscopic structure of shocks in one conservation laws. *Ann. Inst. H. Poincaré Anal. Non Linéaire* **12** 119–153. [MR1326665](MR1326665)
[19] ROST, H. (1981). Nonequilibrium behaviour of a many particle process: Density profile and local equilibria. *Z. Wahrsch. Verw. Gebiete* **58** 41–53. [MR635270](MR635270)
[20] SEPPÄLÄINEN, T. (1998). Coupling the totally asymmetric simple exclusion process with a moving interface. *Markov Process. Related Fields* **4** 593–628. [MR1677061](MR1677061)
[21] SEPPÄLÄINEN, T. (2001). Second class particles as microscopic characteristics in totally asymmetric nearest-neighbor $K$-exclusion processes. *Trans. Amer. Math. Soc.* **353** 4801–4829. [MR1852083](MR1852083)
[22] TALAGRAND, M. (1995). Concentration of measure and isoperimetric inequalities in product spaces. *Inst. Hautes Études Sci. Publ. Math.* **81** 73–205. [MR1361756](MR1361756)
[23] WÜTHRICH, M. V. (2000). Asymptotic behavior of semi-infinite geodesics for maximal increasing subsequences in the plane. In *In and Out of Equilibrium* (V. Sidoravicius, ed.) 205–226. Birkhäuser, Basel.

INSTITUTO DE MATEMÁTICA E ESTATÍSTICA
UNIVERSIDADE DE SÃO PAULO
CAIXA POSTAL 66281
05311-970 SÃO PAULO
BRASIL
E-MAIL: [pablo@ime.usp.br](pablo@ime.usp.br)
E-MAIL: [ordnael@ime.usp.br](ordnael@ime.usp.br)